\documentclass[12pt]{amsart}
\usepackage{latexsym,amsmath,amssymb,amscd}
\usepackage{latexsym,amsmath,amssymb,amscd}
\usepackage{eucal}
\theoremstyle{plain}
\newtheorem{theorem}{Theorem}[section]
\newtheorem{corollary}[theorem]{Corollary}

\newtheorem{proposition}[theorem]{Proposition}

\theoremstyle{definition}
\newtheorem{definition}[theorem]{Definition}

\theoremstyle{remark}

\newcommand{\larray}{\left(\begin{array}{cc}}
\newcommand{\rarray}{\end{array}\right)}

\newcommand{\hattensor}{\hat{\otimes}}

\DeclareMathOperator{\Hom}{Hom}
\DeclareMathOperator{\RH}{H}
\DeclareMathOperator{\HE}{HE}

\DeclareMathOperator{\HP}{HP}

\DeclareMathOperator{\ch}{ch}

\newcommand{\ra}{\rightarrow}

\newcommand{\xr}{\xrightarrow}

\newcommand{\frechet}{Fr\'{e}chet}

\newcommand{\ztwo}{\bbz/2\bbz}
\newcommand{\an}{\text{an}}

\newcommand{\bdsplit}{\begin{displaymath}
\begin{split}}
\newcommand{\edsplit}{\end{split}
\end{displaymath}}

\newcommand{\beqn}{\begin{equation}}
\newcommand{\eeqn}{\end{equation}}
\newcommand{\bsplit}{\begin{split}}
\newcommand{\esplit}{\end{split}}

\newcommand{\bbc}{\mathbb{C}}

\newcommand{\bbz}{\mathbb{Z}}

\newcommand{\calb}{\mathcal{B}}

\newcommand{\calj}{\mathcal{J}}

\newcommand{\calt}{\mathcal{T}}

\newcommand{\fb}{\mathfrak{B}}

\newcommand{\fl}{\mathfrak{L}}

\newcommand{\fs}{\mathfrak{S}}

\usepackage{pb-diagram}
\usepackage{eucal}
\begin{document}
\title[]{Entire cyclic homology of Schatten ideals}
\author{Jacek Brodzki}

\address{School of Mathematical Sciences, University of Southampton, Southampton SO17 1BJ, U.K.}
\email{j.brodzki@soton.ac.uk}
\author{Roger Plymen}
\address{School of Mathematics, University of Manchester,  Manchester 
M13 9PL, U.K.}
\email{roger@maths.man.ac.uk}

\maketitle

Certain cocycles constructed by Connes  are characters of $p$-summable 
Fredholm modules defined over a $\bbc$-algebra $A$ \cite{Connes:PubIHES,Connes:entire}. In this article, we establish some consequences 
of the universal properties which these characters enjoy. Let $\fl^p$ denote the 
$p$-th Schatten ideal, with $1\leq p < \infty$. The characters live in the 
following commutative diagrams. In the case of even Fredholm modules we have
$$
\begin{CD}
\HE^0(\fl^p) @>>> \HE^{0} (A) \\
@A{\simeq}AA @AAA\\
\HP^0(\fl^p) @>>> \HP^{0}(A)
\end{CD}
$$
while in the case of odd Fredholm modules the corresponding diagram is 
$$
\begin{CD}
\HE^0(\fl^p) @>>> \HE^{1} (A) \\
@A{\simeq}AA @AAA\\
\HP^0(\fl^p) @>>> \HP^{1}(A)
\end{CD}
$$
Our main technical result is that the entire cyclic cohomology $\HE^*(\fl^p)$ 
is independent of $p$ and 
isomorphic to 
the entire cyclic cohomology of the ideal of trace class operators $\fl^1$:
$$
\HE^*(\fl^p) \simeq  \HE^*(\fl^1)
$$
We also have 
$$\HE^0(\fl^1) \simeq \bbc
\quad \text{and}\quad \HE^1(\bbc) =0
$$
We establish  companion statements in the entire cyclic homology.

In the context of the above diagram, each character is the image of a
generator of the one-dimensional space $\HP^0(\fl^1)$. Ultimately, each character originates in 
this way.

This observation  relies on a detailed analysis of the entire cyclic cohomology of Schatten ideals $\fl^p$, which 
forms the main technical result of the paper.  This study is made possible by results of 
Meyer \cite{Meyer}, who 
put the entire cyclic cohomology in a more universal bivariant  setting. A key new ingredient 
in his approach is the use of \emph{bornology} instead of topology. He shows 
that the bornology encodes the classical growth condition in the entire cyclic theory 
\cite[Thm. 3.47]{Meyer}.  
At the level of cohomology this 
results in replacing continuous cochains by bounded ones. 

Thanks to our calculations in Section \ref{3}, the method of Cuntz \cite[p. 40]{Cuntz:Exp} can be adapted
to the entire theory. Cuntz's method, which he applied in the case of $kk$-theory
and bivariant periodic cyclic homology,  relies on certain key algebraic properties of those theories 
and on the fact that they both satisfy excision. We show that similar properties can be established 
in the entire case. 

The  $kk$-theory of Cuntz \cite{CDoc}  is universal
among all bivariant functors which are stable,  
half-exact, and invariant with respect to smooth homotopies \cite[p. 48]{Cuntz:Exp}. Thanks to the work of Cuntz \cite[Satz 6.12, p. 173]{CDoc}, we have 
$$kk_0(\bbc, \fl^1)= \bbc.$$
 Thus generators of 
$\HP^0(\fl^1)$ and $\HE^0(\fl^1)$ are the images, via the Chern characters, of generators of
$kk_0(\bbc, \fl^1)$.

\section{Spaces with bornology}

A \emph{bornology} on a set $X$ is a family $\calb$ of subsets of $X$ which is stable 
under the formation of subsets and finite unions \cite{Bou}. Elements of $\calb$ will be 
called \emph{bounded} sets. A \emph{base} of a bornology $\calb$ is a subfamily 
$\calb_0$ of $\calb$ with the property that any element of $\calb$ is contained in some 
element of $\calb_0$. A map $f: (X, \calb) \ra (Y, \calb')$ of bornological 
spaces is \emph{bounded} if and only if for any $B\in \calb$, $f(B) \in \calb'$. 

A bornology on a vector space $E$ is said to be  compatible with the vector space structure
 iff the vector addition $E\times E \ra E$  and multiplication by scalars $\bbc \times E \ra E$ are bounded maps of bornological spaces. 

When the vector space $E$ is equipped with a topology, e.g. $E$ is a locally convex topological 
vector space, then there is a canonical bornology associated with the topology consisting of all 
sets which are absorbed by all neighbourhoods of zero. This is the so called von Neumann 
bornology. However, one can equip the vector space $E$ with a bornology which does not arise 
from the topology of $E$. One consequence of this is that the class of bounded linear maps on $E$ will be different in general from the class of continuous maps. 

The vector space $E$ is called a convex bornological space if it is equipped with a bornology 
whose base consists of convex sets. In this case the base can be assumed to consist of balanced 
convex sets, which will be called \emph{discs}. If $D$ is a bounded disc in $E$ we denote by 
$E_D$ the vector space generated by $D$ and equipped with the seminorm given by the 
gauge of $D$. When $E$ is a 
Hausdorff space, $E_D$ is a normed space.  The spaces $E_D$ form an inductive 
system indexed by the directed family of bounded discs and $E$ is the direct limit of this system. 

We shall say that a disc $D$ is pre-complete iff the space $E_D$ is complete. Thus when $E$ is Hausdorff, $E_D$ is 
a Banach space. We shall say that $E$ is a \emph{complete bornological} space iff 
its bornology admits a base consisting of pre-complete discs (cf. \cite[Def. IV.2.1]{HN}). A bornological 
space $E$ is complete iff it is the inductive limit of an \emph{injective} inductive system of Banach 
spaces \cite[IV.2.3]{HN} \cite[Theorem A.4]{Meyer}. 
Every bornological space $E$ admits a bornological completion, but this operation is less well 
behaved than the usual completion of a uniform space (see \cite[Chapter 4]{HN} \cite[Appendix A]{Meyer} for a full discussion). 

We say that  $A$ is a bornological algebra iff it is equipped with a vector space bornology with respect to which  the product map $A\times A\ra A$ is bounded. A complete bornological algebra 
is one that is complete as a bornological vector space. 

The construction of cyclic type homology theories requires the use of tensor products of 
algebras. In the study of cohomology theories like the entire cyclic cohomology of Banach algebras
it is important to have control over bounded sets in the tensor product of algebras. 
However, as is well known from the work of Grothendieck 
\cite[Probl\'{e}me des Topologies, p.~33]{Gr}, there 
is in general no obvious relation between bounded sets in the tensor product and bounded sets 
in the algebra. 
A possible resolution of  this problem, which works well in some situations, is to define completed 
tensor products with respect to a given bornology rather than topology. 

The bornological tensor product of two bornological spaces $(E_1, \calb_1) $ and $(E_2, \calb_2)$
is by definition the algebraic tensor product $E_1\otimes E_2$ equipped with the bornology whose 
base consists of balanced convex hulls of sets of the form $B_1\otimes B_2$, where $B_1\in 
\calb_1$ and $B_2\in \calb_2$. The completed bornological tensor product $E_1\hattensor
E_2$ is by definition the bornological completion of $E_1\otimes E_2$ with respect this bornology. 

For example, the completed bornological tensor product of two \frechet\ spaces equipped with the 
precompact bornology is isomorphic to  the completed  projective tensor product of the two spaces \cite[Theorem 2.29]{Meyer}. 

When $V_1$ and $V_2$ are \emph{nuclear LF}-spaces regarded as bornological spaces 
with the von Neumann bornology then the completed bornological tensor product $V_1\hattensor V_2$
 is isomorphic  to the inductive tensor product 
$V_1 \bar{\otimes} V_2$ of Grothendieck, see \cite[Cor. 2.30, p.~15]{Meyer}. 

\section{The $X$-complex}

All cyclic type homology theories of an algebra $A$ are defined using a $\bbz/ 2\bbz$-graded 
complex associated with $A$ which, loosely speaking, is constructed using a certain 
deformation of the tensor algebra of $A$. We describe this construction first in the case of 
an algebra $A$ without any topology or bornology. 

We recall the differential graded algebra of differential forms $\Omega A$ associated with $A$. 
By definition, $\Omega A$ is generated by elements of $A$ together with symbols $da$, for $a$ 
in $A$, such that $da$ is linear in $a$ and satisfies the Leibniz rule $d(ab) = (da) b + a db$. 
If the algebra $A$ is unital with unit $1$, it is not assumed that $d(1) = 0$. As a consequence, 
in degree $n$, $\Omega ^n A = A^{\otimes n+1} \oplus A^{\otimes n}$. Elements of $\Omega ^n A$
are linear combinations of differential forms $a_0da_1\dots da_n$ and $da_1\dots da_n$, 
with $a_i$ in $A$. The graded space $\Omega A$ is turned into a differential complex 
by means of two operators
$$
\mathsf{b} = \larray b & 1 - \lambda \\
0 & - b'
\rarray, \qquad \mathsf{B} = \larray 0 & 0 \\ N & 0 \rarray
$$
where the operators $b', b, \lambda, N$ have their usual meaning, c.f. \cite[p.9]{Cuntz:Exp}.

For any algebra $A$ we define the $X$-complex $X(A)$ of $A$ to be the $\ztwo$-graded complex
$$
A \; \mathop{\rightleftarrows}\limits^{\natural d}_{b} \; \Omega^1A_\natural
$$
where $\Omega^1 A_\natural = \Omega^1A/[A , \Omega^1 A] = \Omega^1 A / b(\Omega^2 A)$, 
and $\natural :  \Omega ^1 A \ra \Omega^1 A_\natural $ is the canonical projection map, 
compare \cite[p. 21]{Cuntz:Exp}. 

In order to obtain an interesting homology theory we need to apply the above 
construction to the \emph{non-unital}
 tensor algebra $TA$ of the algebra $A$. While it seems at first sight that 
this will lead to a huge and unwieldy object, it turns out in fact, thanks to the following 
result \cite[Theorem 5.5]{CQ3}, that the resulting complex is a deformation of the mixed complex
$(\Omega A, \mathsf{b}, \mathsf{B})$. A starting point for the proof is a  remark that
 the tensor algebra of any algebra $A$ 
can be identified with the even part of the algebra $\Omega A$ 
 equipped with the Fedosov product 
$$
\omega \circ \eta  = \omega \eta - (-1)^{\deg \omega} d\omega d\eta.
$$
\begin{proposition} (\cite{CQ3}\cite[Thm. 2.29, p. 24]{Cuntz:Exp})
For any algebra $A$ the $X$-complex of the tensor algebra $TA$ of $A$ is isomorphic 
to a complex of the form 
$$
\Omega^{\text{ev}} A \mathop{\rightleftarrows}\limits_{\beta}^{\delta}
\Omega^{\text{odd}} A
$$
where the differentials $\beta$ and $\delta$ can be explicitly determined in terms of 
differentials $\mathsf{b}$ and $\mathsf{B}$.
\end{proposition}

Let us assume that $A$ is a complete bornological algebra. 
Then  $\Omega A$ becomes  a bornological algebra with 
bornology whose base is given by balanced convex hulls of the sets  $S dS \dots dS$ and $dS \dots dS$, where $S$ is an element of the bornology on $A$. We denote by $\Omega_{\text{an}}A$ the completion 
of $\Omega A$ with respect to this bornology. This is a $\ztwo$-graded complex
with the same differentials $\mathsf{b}$ and $\mathsf{B}$, which are now bounded maps. 
The even part of the algebra $\Omega_{{\an}}A $, equipped with the Fedosov product, 
is by definition the analytic tensor algebra $\calt A$ of $A$. This algebra fits into the 
algebra extension 
$$
0 \ra \calj  A\ra \calt A \ra A\ra 0
$$
The $X$-complex of the tensor algebra $\calt A$ is defined in the same way 
as before.

\begin{definition} The bivariant entire
cyclic homology of a pair of bornological algebras $A$ and $B$ is by definition 
$$
\HE_*(A, B) = H_*(\Hom(X(\calt A), X(\calt B))
$$
where $\Hom(X(\calt A), X(\calt B))$ denotes the complex  of bounded linear maps from 
$X (\calt A)$ to $X(\calt B)$. This complex is  equipped with the differential $[\partial , \phi] = \partial \circ 
\phi - (-1)^{\deg{\phi}}\phi\circ\partial$, where $\partial = \mathsf{b}+ \mathsf{B}$
 \cite[p.~57]{Cuntz:Exp}\cite[p.~37]{Meyer}.
\end{definition}
 This construction  is due to Meyer, 
who proves that this bivariant cyclic homology satisfies excision in both variables \cite{Meyer}
\cite[Thm. 5.4]{Cuntz:Exp}. More precisely, we have the following. 
\begin{theorem}\label{excision}
Let $0\ra S\ra P\ra Q\ra 0$ be an extension of complete bornological algebras which admits 
a bounded linear section. Assume further that $A$ is a complete bornological algebra. 
Then we have the following natural exact sequences of length six. 
\begin{equation}\label{firstsequence}
\begin{CD}
\HE_0(A, S) @>>> \HE_0(A, P)@>>> \HE_0(A, Q)\\
@A\mathsf{d}_1AA @. @VV\mathsf{d}_1V\\
\HE_1(A,Q) @<<< \HE_1(A,P) @<<< \HE_1(A,S)
\end{CD}
\end{equation}
\begin{equation}\label{secondsequence}
\begin{CD}
\HE_0(S,A) @<<< \HE_0(P,A) @<<< \HE_0(Q,A) \\
@V\mathsf{d}_2VV @.@AA\mathsf{d}_2A\\
\HE_1(Q,A) @>>> \HE_1(P, A) @>>> \HE_1(S,A)
\end{CD}
\end{equation}
\end{theorem}
Moreover, $\HE$ 
is invariant  with respect to  differentiable  homotopies whose  first derivative is integrable. Meyer also proves that 
when $A$ is a Banach algebra then $\HE_*(A, \bbc)$ is the same as the entire cyclic 
cohomology $\HE^*(A)$ of $A$ as defined by Connes \cite[4.1]{Meyer:Excision}

A very important property of $\HE$ is the existence of the composition product, which is 
defined as in the case of bivariant periodic cyclic homology by composition of linear chain 
maps. For any three  bornological algebras $A_1$, $A_2$ and $A_3$ there 
is a bilinear product: 
$$
\HE_i(A_1,A_2) \times \HE_j(A_2, A_3) \ra \HE_{i+j}(A_1, A_3)
$$
given by $f\cdot g = g\circ f$. A detailed study of this product will provide formulae for the connecting
homomorphisms in the exact sequences of Theorem \ref{excision}.

\section{Entire cyclic cohomology of Schatten ideals}\label{3}

The aim of this section is the following 
\begin{theorem}\label{3.1}
Let $\fl^p$ and $\fl^q$ be two Schatten ideals, $ 1 \leq p < q$. Then  the inclusion $\fl^p \ra \fl^q $ 
induces an invertible element in $\HE_*(\fl^p, \fl^q )$. Consequently, 
the entire cyclic homology and cohomology of the two algebras are isomorphic: 
$$
\HE^*(\fl^p) = \HE^*(\fl^q); \qquad 
\HE_*(\fl^p ) = \HE_*(\fl^q)
$$
\end{theorem}

In the context of algebraic periodic cyclic homology \cite{Cuntz:Exp} and in $kk$-theory 
this result was first proved by Cuntz \cite{CDoc}. The proof outlined here follows the same strategy, 
which relies on algebraic features of bivariant cohomology theories. To make sure that 
this translation works, we need to prove a number of technical results which provide the necessary formal 
properties of the bivariant entire cyclic homology. 

Let $E$ denote the following extension of complete bornological algebras. 
$$
E: \qquad 0 \ra S \xr{i} P \xr{p} Q\ra 0
$$
We shall assume that this sequence is split, i.e. there exists a bounded linear map $s : Q \ra P$  which 
is a right inverse  for the projection $p$. 

The excision property of the bivariant cyclic homology $\HE_*$ implies that there are the following 
two exact sequences 
\begin{displaymath}
\begin{split}
\ra & \HE_*(P, S) \ra \HE_*(S,S) \xr{\delta_1} \HE_{*+1}(Q, S) \ra \\
\ra & \HE_*(Q, P) \ra \HE_*(Q,Q) \xr{\delta_2} \HE_{*+1}(Q, S) \ra 
\end{split}
\end{displaymath}
Denote by $1_Q$ the class in $\HE_0(Q,Q)$ induced by the identity map on the algebra $Q$ 
and, similarly, $1_S$ will denote the class of the identity map on $S$. The following lemma is a 
translation of  a result of Kassel \cite[Lemme 2.2]{Kassel} in the case of his
bivariant cyclic cohomology. An analogous result  has
been proved in the case of bivariant periodic cyclic cohomology by Cuntz and Quillen 
\cite{CQ:Excision}\cite[Prop. 2.51, p.33]{Cuntz:Exp}. 
\begin{proposition}\label{3.2}
If $\delta_1$ and $\delta _2$ denote the connecting homomorphisms in the preceding two long exact
 sequences, then
$$
\delta_1(1_S) = - \delta_2(1_Q) \in \HE_1(Q,S)
$$
\end{proposition}
\begin{proof}
To simplify notation, for any two $\bbz/2\bbz$-graded complexes $C$ and $D$, 
we shall write $\RH_*(C,D)$ for the homology $\RH_*(\Hom(C,D))$. 

The projection $p$ induces a map of complexes $X (\calt P) \ra X(\calt Q)$. Let us denote by 
$X(P,Q)$ the kernel of this map, so that we have an exact sequence of $\ztwo$ graded complexes
$$
(\alpha): \quad 0 \ra X(P,Q) \ra X(\calt P) \ra X(\calt Q) \ra 0.
$$
By the universal properties of the non-unital tensor algebra and the $X$-complex, 
a linear splitting $s$ of the sequence $E$ induces a linear splitting of the sequence $(\alpha)$ 
\cite[3.3.2]{Meyer} which will also be denoted $s$. This splitting in 
turn determines the odd degree map $[\partial, s] = \partial \circ s - s \circ \partial$ 
with the properties $p[\partial , s] = 0$
and $[\partial, [\partial, s]]=0$. Thus $[\partial, s]$ determines an element in $\RH_1(X(\calt Q), 
X(P, Q))$ which will be denoted  $\gamma$. This class does not depend on the 
choice of the linear section $s$. Indeed, any two $\bbc$-splittings can be connected by a linear 
path $(1-t)s + ts'$. Now:
$$
[\partial, (1-t) s + ts'] = (1-t) [\partial, s] +t [\partial , s']
$$
gives a (differentiable) homotopy between the corresponding cycles. 

Furthermore, there is the induced sequence of $\ztwo$-graded  complexes
\begin{displaymath}
\begin{split}
  0 \ra & \Hom (X (\calt Q), X(P,  Q)) \ra \Hom(X(\calt P), X(P, Q)) \\
& \ra \Hom(X(P,Q), X(P, Q)) \ra 0
\end{split}
\end{displaymath}
If we now apply the homology functor
 we obtain an exact sequence of length six. The two connecting homomorphisms 
 $$
 \RH_j(X(P,Q), X(P,Q)) \xr{\gamma \cdot} \RH_{j+1}(X(\calt Q), X(P,Q)),
 $$
for $j=0,1$,  in the resulting  homology sequence are both given by multiplication by 
$\gamma$. 

The inclusion map $i: X(\calt S) \ra X (\calt P)$ satisfies $pi = 0$ and so can be regarded as 
a  map $i : X(\calt S) \ra X(P,Q)$. Since this map is induced from an algebra 
homomorphism, it is a  $0$-cycle and so creates an element $i \in 
\RH_0(X(\calt S), X(P, Q))$.   A key step in the proof of excision in $\HE$ (Theorem 
\ref{excision}) is the fact 
that $i$ is  invertible \cite{Meyer}, so that there exists $i^{-1}\in \RH_0(X(P,Q), X(\calt S))$. 

For any $\bbz/2\bbz$-graded complex $C$ the composition product gives a map
$$
\RH_j(X(\calt S), X(P,Q))\otimes \RH_k(X(P,Q), C) \ra \RH_{j+k}(X(\calt S), C). 
$$
Thus taking the product on the left by the invertible element $i$ of degree 0 establishes an 
isomorphism
$$
i\cdot : \RH_j(X(P, Q), C) \xr{\simeq} \RH_j(X(\calt S), C).
$$
If we now recall that $\HE_*(S,S) = \RH_*(X(\calt S), X(\calt S))$ and use the triple 
product 
\begin{displaymath}
\begin{split}
\RH_0(X(P,Q), X(\calt S))\otimes & \RH_j(X(\calt S), X(\calt S)) 
\otimes \RH_0(X(\calt S), X(P,Q))\\  & \longrightarrow \RH_j(X(P,Q), X(P,Q))
\end{split}\end{displaymath}
we deduce that there is an isomorphism
$$
\HE_j(S,S) \simeq \RH_j(X(P, Q), X(P,Q)).
$$
For $\phi \in \HE_*(S,S)$  this is given explicitly by $\phi \mapsto i^{-1} \cdot \phi \cdot i$. Taking 
the left product with $\gamma$ produces a map 
$$
\HE_j(S,S) \ra \RH_{j+1}(X(\calt Q), X(P,Q)), \quad \phi\mapsto \gamma\cdot i ^{-1}\cdot \phi 
\cdot i.
$$
If we now take the product on the right with $i^{-1}\in \RH_0(X(P,Q), X(\calt S))$ and use the 
identification $\RH_*(X(\calt Q), X(\calt S)) = \HE_*(Q,S)$ we 
conclude that   the connecting homomorphism $\HE_j(S,S)\ra
\HE_{j+1}(Q,S)$ may be described by the fomula
$$
\phi\mapsto \gamma\cdot i^{-1}\cdot \phi, 
$$
for any $\phi\in \HE_j(S,S)$. In particular, when $\phi = 1_S$ we obtain
$$
\delta_1 (1_S)  = \gamma\cdot i^{-1}. 
$$

The exact sequence  $(\alpha)$ also  leads to the following exact sequence: 
\begin{displaymath}
\begin{split}
0 \ra & \Hom(X(\calt Q), X(P, Q)) \ra \Hom(X(\calt Q), X(\calt P)) \\
& \ra \Hom(X(\calt Q), X(\calt Q)) \ra 0
\end{split}
\end{displaymath}
which in turn induces  a corresponding exact sequence of bivariant 
homology groups of length six.  

The existence of a bounded linear splitting $s$ of the sequence $(\alpha)$ implies that
$X(\calt P)$ splits as a direct sum of $\bbz/2\bbz$-graded vector spaces: 
$$
X(\calt P) = X(\calt Q) \oplus X(P, Q)
$$
It then follows that there is a complementary splitting of the sequence $(\alpha)$ given by 
$s' = 1-s$. It is now not difficult to see that the two connecting homomorphisms
$$
\RH_j(X(\calt Q), X(\calt Q)) \ra \RH_{j+1} ( X(\calt Q), X(P,Q))
$$
for $j=0,1$ are given by left multiplication by the element $[\partial , s'] =  - [\partial , s] = -\gamma$. 

We can find an explicit formula for the connecting 
homomorphism
$$
\delta_2 : \HE_j(Q,Q) \ra \HE_{j+1}(Q,S)
$$
as follows. 
We  multiply on the 
left by  
 $- \gamma$ to construct a map
$$
\HE_j(Q,Q) := \RH_j(X(\calt Q), X(\calt Q)) \xr{-\gamma\cdot} \RH_{j+1} ( X(\calt Q), X(P,Q)). 
$$
Multiplication on the right by 
$i ^{-1} \in\RH_0(X(P,Q), X(\calt S))$ gives a map
$$
\RH_{j+1} ( X(\calt Q), X(P,Q)) \xr{\cdot i^{-1}} \RH_{j+1}(X(\calt Q), X(\calt S)) = \HE_{j+1}(Q,S). 
$$
Thus the formula for the connecting homomorphism $\delta_2$ is 
$$
\delta_2(\psi) = - \gamma \cdot \psi \cdot i^{-1},  
$$
for any $ \psi\in \HE_j(Q,Q)$. 
In particular, when $\psi = 1_Q$ we have 
$$
\delta_2(1_Q) = - \gamma \cdot i^{-1} = - \delta_1(1_S)
$$
\end{proof}

This result can be extended to provide formulae for 
connecting homomorphisms in exact sequences  of Theorem \ref{excision}. Both 
of the excision exact sequences are natural. In the case of the 
sequence (\ref{firstsequence}) this means that there exists a commutative 
diagram
$$
\begin{CD}
\HE_j(A,Q) \times \HE_0(Q,Q) @>m>> \HE_j(A,Q)\\
@V{1\otimes\delta_2}VV @VV\mathsf{d}_1V\\
\HE_j(A,Q) \times \HE_1(Q,S)@>m>> \HE_{j+1}(A,S)
\end{CD}
$$
where $m$ denotes the product map and $\mathsf{d}_1$ denotes the connecting 
homomorphism in the diagram (\ref{firstsequence}) of Theorem \ref{excision}, for $j=0,1$. 

Taking into account the usual sign convention we have that 
$$
\mathsf{d}_1(\phi\cdot \psi) = m(1\otimes \delta_2)(\phi\otimes \psi)
 = m((-1)^{\deg(\phi)}(\phi\otimes \delta_2(\psi)) = 
 (-1)^{\deg (\phi)}\phi\cdot\delta_2(\psi)
 $$
 for $\phi\in \HE_j(A,Q)$ and $\psi\in \HE_0(Q,Q)$. Hence
 $$
 \mathsf{d}_1(\phi) = \mathsf{d}_1(\phi\cdot 1_Q)
 = (-1)^{\deg(\phi)} \phi\cdot\delta_2(1_Q)
 $$
 Similarly, we obtain a formula for the connecting homomorphism $\mathsf{d}_2$
  in the exact sequence 
 (\ref{secondsequence}). In this case the naturality of this sequence implies that there 
 exists the following commutative diagram
 $$
 \begin{CD}
 \HE_0(S,S) \times \HE_j(S,A) @>m>>\HE_j(S,A)\\
 @V{\delta_1}\otimes 1VV @VV{\mathsf{d}_2}V\\
 \HE_1(Q,S) \times \HE_{j}(S,A) @>m>> \HE_{j+1}(Q,A)
 \end{CD}
 $$
 Hence, for $\phi\in \HE_0(S,S)$ and $\phi\in \HE_j(S,A)$ we have
 $$
 \mathsf{d}_2(\phi\cdot\psi) = m(\delta_1\otimes 1)(\phi\otimes \psi)
 = \delta_1(\phi)\cdot \psi
 $$
 Thus 
 $$
 \mathsf{d}_2(\psi) = \mathsf{d}_2(1_S\cdot\psi) = \delta_1(1_S)\cdot \psi.
 $$
 
In summary, we have obtained the proof of the following proposition, which extends an analogous result 
of Kassel \cite[Thm  2.1, Lemme 2.2]{Kassel} (see also  \cite[Thm 5.5]{CQ:Excision}). 
\begin{proposition}\label{3.3}
Let us denote by $\ch(E)$ the class $ -\delta_1(1_S) = \delta_2(1_Q)$ of the extension $E$. 
Then the connecting homomorphism $\mathsf{d}_1$ in the exact sequence (\ref{firstsequence})
sends $\phi \in \HE_j(A,Q)$ to $(-1)^{\deg(\phi)}\phi \cdot \ch(E)\in \HE_{j+1}(A,S)$. The connecting homomorphism 
$\mathsf{d}_2$ in the sequence (\ref{secondsequence}) sends  $\psi \in \HE_j(S,A)$ 
to $\ch(E)\cdot \psi \in\HE_{j+1}(Q,A)$. 
\end{proposition}
This implies, as in \cite{Kassel}\cite{CQ:Excision}, the following. 
\begin{corollary} \label{3.4}
If the algebra $P$ in the extension $E$ is $\HE$-equivalent to $0$, which 
means that $\HE_*(A,P)= \HE_*(P,A)=0$ for any bornological algebra $A$, then $\ch(E)$ is an 
invertible element in $\HE_1(Q,S)$. 
\end{corollary}
\begin{proof}
Let us put $A= S$ in the sequence (\ref{firstsequence}) and then 
$A=Q$ in the sequence (\ref{secondsequence}) of Theorem 
\ref{excision}.  Since the terms containing the 
algebra $P$ are zero, we see that the connecting homomorphisms $\mathsf{d}_1$
and $\mathsf{d}_2$ are now isomorphisms. In particular, there 
exists $\eta_1 \in \HE_1(S,Q)$ such that $\mathsf{d}_1(\eta_1) = 1_S \in \HE_0(S,S)$. 
Similarly, there exists $\eta_2\in \HE_1(S,Q)$ such that 
$\mathsf{d}_2(\eta_2) = 1_Q \in \HE_0(Q,Q)$. 
But we have just established that 
$$
\mathsf{d}_1(\eta_1) = \eta_1 \cdot \ch(E) = 1_S
$$
and that 
$$
\mathsf{d}_2(\eta_2) = \ch(E) \cdot\eta_2 = 1_Q
$$
This implies that $\eta_1 = \eta_2$. Indeed, 
$$
\eta_1= \eta_1\cdot 1_Q= \eta_1\cdot \ch(E) \cdot \eta_2 = 1_S\cdot\eta_2=\eta_2.
$$
Thus $\eta= \eta_1= \eta_2 \in \HE_1(S,Q)$ is the inverse of $\ch(E) \in 
\HE_1(Q,S)$. 
\end{proof}

The rest of the proof of the theorem follows an argument of Cuntz
 \cite[Satz 6.12]{CDoc}\cite{Cuntz:Exp}, who used it to prove an analogous result in his $kk$ theory and periodic cyclic homology.
\begin{proposition}
Let us assume that for two complete bornological algebras $A$ and $B$ there are maps 
$$
\begin{array}{cl}
\alpha : & B \hookrightarrow A \\
\beta: & A\otimes A \ra B
\end{array}
$$
such that the composition $\alpha \circ \beta$  identical to the product map on $A$, whereas 
$\beta \circ \alpha\otimes \alpha$ is the product on $B$. Then 
the element $[\alpha]$ of $\HE_0(B,A)$ is invertible. This implies that 
$\HE^*(A) \simeq \HE^*(B)$ and $\HE_*(B) \simeq \HE_*(A)$. 
\end{proposition}

\begin{proof}
We present here a 
more explicit version of Cuntz's  argument, which is adapted to the context of entire cyclic homology. 

We equip the Fr\'{e}chet algebra $C^\infty([0,1]) $ of smooth functions with the von Neumann
bornology; in the present case this bornology coincides with the pre-compact bornology. 

If $A$ is a complete bornological algebra, define
$$
A[0,1] = C^\infty([0,1]) \hattensor A
$$
We denote by $A(0,1]$ the algebra of smooth functions from the \emph{closed} interval 
$[0,1]$ to $A$ which vanish at zero; we use the notation $A(0,1)$, $A[0,1)$ to denote the algebras
of smooth functions from the interval $[0,1]$ to $A$ that vanish at both ends of the interval or just 
at $1$. The algebra $A[0,1)$ is contractible to zero: the family of maps $\phi_t$ that send 
a function $f$ to $\phi_t(f)(x) = f((1-t)(x))$ forms a homotopy between the identity map and evaluation 
at $1$ (which is the same as the zero map). 

There is the  following suspension extension: 
$$
\fs(A): 0 \ra A(0,1) \ra A[0,1) \ra A\ra 0
$$
where the map on the right is given by evaluation at $0$. Since the algebra $A[0,1)$ 
is contractible it is $\HE$-equivalent to zero. We can therefore use Corollary \ref{3.4} to deduce that 
the class $\ch{\fs(A)}= -\delta_1(1_{A(0,1)}) = \delta_2(1_A)\in \HE_1(A, A(0,1))$ 
is invertible, for any complete 
bornological algebra $A$. 

Let $A$ and $B$ be complete bornological algebras as in the statement of the Proposition. We denote 
by $\calb$ the complete bornological 
 algebra generated by the algebra $B(0,1)$ together with the algebra $At= 
\{f_a\mid a\in A\}$ 
consisting of functions $f_a : [0,1] \ra A $ which for a fixed $a\in A$ send $ t \mapsto ta$.  
As a vector space, $\calb$ is the direct sum of the two algebras. The product on $\calb$ is defined 
using the pointwise product on $B(0,1)$ together with the following two operations. The product 
of a function $f \in B(0,1)$ by an element  $g_a\in At$ is given by $\mu(\alpha(f) \otimes g_a)$. Finally, 
the product of two fucntions $f_a$ and $f_b$ is the function $g(a,b) + f_{\alpha\mu(a\otimes b)}$
where 
$$
g(a,b) (t) = \mu(a \otimes b) (t^2 - t). 
$$
With these definitions we have the following extension of complete bornological algebras
$$
0 \ra B(0,1) \ra \calb \ra A \ra 0
$$
which admits a bounded linear splitting. Proposition \ref{3.2} implies that this extension creates an
element $u \in \HE_1(A, B(0,1))$. 

The homomorphism $\alpha: B\ra A$ gives rise to an element $[\alpha] \in 
\HE_0(B,A)$. We need to show that it is invertible. For this we construct first the following
diagram:
$$
\begin{CD}
0 @>>> B(0,1) @>>> B[0,1) @>>> B @>>> 0\\
@. @VVV @VV{id + \alpha}V @VV{\alpha}V\\
0 @>>> B(0,1) @>>>\calb @>>> A @>>> 0
\end{CD}
$$
Using the first of the two excision sequences we obtain the following commutative diagram:
$$
\begin{CD}
\HE_0(B(0,1]) @>>> \HE_0(B) @>{\cdot\ch(\fs(B))}>> \HE_1(B(0,1)) @>>>\\
@VVV @V{\alpha}VV @| \\
\HE_0(\calb) @>>> \HE_0(A) @>{\cdot u}>>\HE_1(B(0,1)) @>>>
\end{CD}
$$
Here the two connecting homomorphisms on the right are in accordance with Proposition 
\ref{3.3} whereas the vertical map in the middle is given by taking the product on the right with 
$\alpha \in \HE_0(B,A)$. Since the diagram commutes we see that for any $\phi 
\in \HE_0(B)$
$$
\phi\cdot \alpha\cdot u = \phi\cdot \ch(\fs(B))$$
Given that $\ch(\fs)$ is invertible we find that $\alpha \cdot u \cdot \ch(\fs(B))^{-1} = 1\in \HE_0(B,B)$.

We now employ the following commutative diagram:
$$
\begin{CD}
0 @>>> B(0,1) @>>>\calb @>>> A @>>> 0\\
@. @V{\alpha}VV @VV{\alpha}V @| \\
0 @>>> A(0,1) @>>> A [0,1) @>>> A @>>> 0
\end{CD}
$$
where the vertical map on the left is the obvious extension of $\alpha$ to functions. 
Using excision again, this translates to the following commutative diagram of homology groups: 
$$
\begin{CD}
\HE_0(\calb) @>>> \HE_0(A) @>{\cdot u}>>\HE_1(B(0,1)) @>>>\\
@V{\alpha}VV @| @VV{\alpha}V \\
\HE_0(A[0,1)) @>>> \HE_0(A) @>{\cdot \ch(\fs(A))}>> \HE_1(A(0,1)) @>>>
\end{CD}
$$
If we take into account the isomorphism of homology groups $\HE_1(B(0,1)) = \HE_0(B)$ 
provided by the extension extension (and similarly in the case of the bottom row) we 
obtain the commutative diagram: 
$$
\begin{CD}
\HE_0(A) @>{\cdot u}>>\HE_1(B(0,1)) @>{\cdot \ch(\fs(B))^{-1}}>> \HE_0(B) \\
@| @VV{\alpha}V @VV{\alpha}V \\
\HE_0(A) @>{\cdot \ch{\fs(A)}}>> \HE_1(A(0,1)) @>{\cdot \ch(\fs(A))^{-1}}>> \HE_0(A) 
\end{CD}
$$
It is now clear that $u\cdot \ch(\fs(B))^{-1} \cdot \alpha = 1 \in \HE_0(A,A)$. Thus 
$\alpha \in \HE_0(B,A)$ is invertible, with inverse $u\cdot \ch(B(0,1))^{-1} \in 
\HE_0(A,B)$. 
\end{proof}

To finish the proof of the theorem we 
let  $B= \fl^p$ and $A = \fl^q$, where $p\leq q \leq 2p$. The map $\alpha$ of the previous statements 
is obtained from the continuous inclusion $\fl^p \ra \fl^q$ and the map $\beta $ from the 
multiplication map $\fl^q \hattensor \fl^q \ra \fl^p$. This completes the proof of Theorem \ref{3.1}. 

\begin{corollary}
 Let $\fb$ be a complete
bornological algebra. Then for any $1\leq p < q$ the inclusion $\fl^p\hattensor \fb \ra 
\fl^q \hattensor \fb$ induces an invertible element in bivariant entire cyclic cohomology
$\HE^0(\fl^p \hattensor \fb, \fl^q\hattensor \fb)$. Thus the entire cyclic homology and cohomology 
of the algebras  $ \fl^p \hattensor \fb$ and $\fl^q \hattensor \fb$ are isomorphic. 
\end{corollary}

\section{Hochschild homology of $\fl^1$} 

In the previous section we have established $\HE^*$-equivalence of the Schatten ideals $\fl^p$
for all $p\geq 1$. In this section we prove that the ideal of trace class operators 
$\fl^1$ is $\HE^*$-equivalent to $\bbc$.  Although this was proved in \cite{Meyer} in the 
context of the bivariant theory, we present here a simple direct proof which relies on 
results about Hochschild homology of $\fl^1$.

Let $E$ and $F$ be two Banach spaces in duality relative to a non-degenerate bilinear form
$\langle - , - \rangle: E \times F \ra \bbc$. Then the tensor product of these spaces can 
be turned into an algebra with the multiplication defined by 
$$
(x\otimes y) (x'\otimes y') = \langle x', y'\rangle x \otimes y'
$$
In the case where $E$ is a Hilbert space $H$ and $F$ is its continuous dual $H^*$, we have 
that $H\hattensor H^*= N(H)$, where $N(H)$ is  the algebra of nuclear operators on $H$. 
When $H$ is separable, the algebra of nuclear 
operators is isomorphic to the algebra of trace class operators $\fl^1$. 

Furthermore, Helemskii proves in \cite[Ch.~IV]{Helem} that the algebra of nuclear operators $N(H)$, hence the 
algebra of trace class operators $\fl^1$ is biprojective, which  means that 
it is a projective bimodule over itself. It is also proved in  \cite[Theorem V.2.28]{Helem} that for 
a biprojective Banach  algebra $A$ we have $H^3(A, X) = 0$ for any Banach bimodule $X$. 
It then follows   \cite{K} that there exists a connection 
that provides a uniformly bounded contracting homotopy of the Hochschild complex for $A$ (with coefficients in 
$A$). This implies by the perturbation mapping lemma that
the canonical inclusion $\HP^*(A) \ra \HE^*(A)$ is an isomorphism \cite{K}. On the left-hand side 
of this map,
we regard $A$ as a \emph{topological} algebra and 
define $\HP^*(A)$  via the projective tensor product. 

To summarise, this sequence of arguments shows that $\HP^*(\fl^1) \simeq \HE^*(\fl^1)$. 
Finally, Cuntz proves in \cite[Prop. 17.3]{Cuntz:Exp} that $\HP^*(\fl^1) =\HP^*(\bbc)$ and $\HP_*(\fl^1) 
=\HP_*(\bbc)$. 
We summarise these results as follows. 
\begin{proposition}
The algebra $\fl^1$ is $\HE$-equivalent to $\bbc$. It follows from Theorem \ref{3.1} that the same is 
true for any $\fl^p$, $p\geq 1$. 
\end{proposition}

\section{Canonical classes associated with $p$-summable Fredholm modules} 

In this section we use our calculations to put some known results concerning characters of 
Fredholm modules in a new context. Let $A$ be an involutive algebra over $\bbc$. 

We begin with the odd case. We recall that an odd $p$-summable Fredholm module over 
$A$ is given by the data
$(H, \pi, F)$, where $\pi: A \ra \fl(H)$ is a representation the  algebra $A$ on a Hilbert space $H$, and $F$ is a self-adjoint involution which 
commutes with $\pi$ modulo $\fl^p$ \cite[p. 208]{Connes:Book}. 

Let $P$ be the
corresponding  spectral projection onto the $+1$ eigenspace. Let $\sigma: A \ra \fl(H)$ be a linear
map defined by $\sigma(a) = P\pi(a) P$ for all $a\in A$, where $\pi$ is the representation of $A$ 
as bounded operators on Hilbert space as required by the structure of a Fredholm module. The goal 
of this section is to construct canonical classes in the periodic and entire cohomology of the 
algebra $A$. Our construction relies on an idea of Cuntz and Quillen \cite{CQCRAS}, and 
follows the method outlined by Cuntz in \cite{CDoc,Cuntz:Exp}. 

Let $A' = \fl^p +\sigma(A)$. This is a subalgebra of $\fl (H)$. 
The Schatten ideal $\fl^p$ is then an  ideal in the algebra $A'$ and we have the following short
exact sequence of algebras, which is $\bbc$-split:
\begin{equation}\label{extension1}
0 \ra \fl^p \ra A' \ra A'/\fl^p \ra 0
\end{equation}
The linear map $\sigma$ can be viewed as a map $\sigma : A \ra A'$, which gives rise to 
an algebra homomorpism $\sigma: TA \ra A'$ which has the important property that it 
sends the canonical ideal $IA \subset TA$ to the Schatten class $\fl ^p$. In other words, we have the following commutative diagram of short exact sequences: 
$$ 
\begin{CD}
0 @>>> \fl^p @>>>  A' @>>> A'/\fl^p @>>> 0 \\
@. @AAA @AAA @AAA \\
0 @>>> JA @>>> TA @>>> A @>>> 0 
\end{CD}$$
The  algebra homomorphism $JA \ra \fl ^p$ gives rise to an element in $\HP_0(JA, \fl^p)= 
\HP_1(A,\fl^p)$ 
and so to a map in cohomology $\HP^*(\fl ^p) \ra \HP^*(JA) = \HP^{*+1}(A)$. 
Given that $\HP^0(\fl^p) = \bbc$ and $\HP^1(\fl^p) = 0$, an 
odd Fredholm module determines a canonical element in $\HP^1(A)$. This is the character of 
a $p$-summable Fredholm module that was first constructed by Connes in \cite{Connes:PubIHES}. 

Our discussion of bornological algebras allows us to extend this idea to entire cyclic cohomology. 
Let $A$ be a complete bornological algebra. Applying the same reasoning as above to the 
canonical extension of complete bornological algebras
$$
0 \ra \calj A \ra \calt A \ra A \ra 0
$$
produces an 
element of $\HE_1(\calj A, \fl^p)$ and so a map $\HE^*(\fl^p) \ra \HE^*(\calj A) = \HE^{*+1}(A)$. 
Again, since we have that $\HE^0(\fl^p) = \bbc$ and $\HE^1(\fl^p) = 0$, an odd $p$-summable 
Fredholm module determines a canonical class in $\HE^1(A)$. Furthermore, because of 
 of $\HP$ and $\HE$ equivalence of Schatten ideals, these classes are independent of 
 $1 \leq p < \infty$. 
 
We remark that the canonical classes so constructed are compatible, in the sense that 
we have the following commutative diagram
$$
\begin{CD}
\HE^0(\fl^p) @>>> \HE^1 (A) \\
@A{\simeq}AA @AAA\\
\HP^0(\fl^p) @>>> \HP^1(A)
\end{CD}
$$
It is not difficult to work out the well-known explicit formulae for these characters; these were first derived by Connes in \cite{Connes:PubIHES}, compare \cite[Ch. 19]{Cuntz:Exp}. 

An \emph{even} $p$-summable Fredholm module over an involutive $\bbc$ algebra $A$ is given by the data $(H, \pi, F, \gamma)$, where 
$\gamma$ is a self-adjoint involution on the Hilbert space $H$ (this Hilbert space is thus 
$\bbz/2\bbz$-graded) and the representation $\pi : A \ra \fl (H)$ commutes with this involution (and
so $A$ is represented  by even operators with respect to the grading). $F$ and $\gamma$ anticommute
and for each $a\in A$,  $[F, \pi(a)]\in \fl^p$. 

In this situation we have a different algebra extension: 
\begin{equation}\label{extension}
0\ra \fl^p \ra A_\gamma \ra A_\gamma/\fl^p \ra 0
\end{equation}
 where the algebra $A_\gamma$ is generated by $\fl^p$ and $\pi(A)$. This sequence has two 
linear splittings: $\pi$ and $\pi^F(a) = F\pi(a)F$, see \cite[Ch. 19]{Cuntz:Exp}. In the context on 
periodic cyclic cohomology this extension leads to the well known canonical character of Connes:
$$
\HP^0(\fl^p) \ra \HP^0(A)
$$
which again is independent of $p$. This construction carries over to the case when the algebra 
$A$ is a complete bornological algebra; in particular the algebra extension (\ref{extension}) becomes 
an extension of complete bornological algebras. This extension gives rise to a map
$\HE^0(\fl^p) \ra \HE^0(A)$ which represents the character of an even Fredholm module. 
The two constructions, one in the context of periodic cyclic cohomology and the other for the 
entire cyclic cohomology, are compatible in the sense that there exists the following commutative 
diagram: 
$$
\begin{CD}
\HE^0(\fl^p)  @>>> \HE^0(A) \\
@A{\simeq}AA @AAA\\
\HP^0(\fl^p) @>>> \HP^0(A)
\end{CD}
$$
These remarks may be summarised as follows. Let $A$ be a complete involutive bornological 
algebra over $\bbc$. Let $\alpha$ be an odd $p$-summable Fredholm module over $A$, and
let $\alpha_\gamma$ be an even $p$-summable Fredholm module over $A$. Let 
$\ch_E(\alpha)\in \HE^1(A)$ be the class  in $\HE^1(A)$ determined by $\alpha$, and 
$\ch_P(\alpha)\in \HP^1(A)$ be the class in $\HP^1(A)$ determined by $\alpha$. 
Similarly, we denote by $\ch_E(\alpha_\gamma)\in\HE^0(A)$ and 
$\ch_P(\alpha_\gamma)\in\HP^0(A)$
the classes in $\HE^0(A)$ and $\HP^0(A)$ determined by $\alpha_\gamma$. 
\begin{theorem}
Under the canonical inclusion 
$$\HP^*(A) \ra \HE^*(A),$$
$\ch_E(\alpha)$ is the image of  $\ch_P(\alpha)$. 
In the even case,  $\ch_E(\alpha_\gamma)$ is the image of
$\ch_P(\alpha_\gamma)$. 
\end{theorem}
\begin{corollary}
The class $\ch_E(\alpha)\in \HE^1(A)$ can be represented by a periodic cyclic cocycle. 
The class $\ch_E(\alpha_\gamma) \in \HE^0(A)$ can be represented by a periodic cyclic cocycle. 
\end{corollary}

\end{document}